\documentclass{svproc}
\usepackage{url}
\usepackage{color}
\usepackage{bm}
\usepackage{ulem}
\usepackage{comment}

\usepackage{siunitx}
\usepackage{amsmath,amssymb}
\usepackage{hyperref}
\usepackage[capitalise]{cleveref}

\newcommand{\numgru}{\num[group-minimum-digits=3]}
\usepackage{textcomp}
\usepackage{graphicx}
\usepackage{multicol}
\usepackage{multirow} 
\usepackage{floatrow}
\newfloatcommand{capbtabbox}{table}[][\FBwidth]

\usepackage{blindtext}
\definecolor{darkgreen}{RGB}{0,120,0}

\newcommand{\diff}{\,\mathrm{d}}

\usepackage{ifthen}
\newcommand{\mycomment}[1]{%
\ifthenelse{\isodd{\value{page}}}{%
\normalmarginpar%
\marginpar{\tiny {#1}}%
}{%
\reversemarginpar%
\marginpar{\tiny {#1}}%
}}%

\begin{document}
\mainmatter              % start of a contribution
%
%\title{GDSW coarse spaces for finite element simulations of laser beam welding processes}
\title{Monolithic overlapping Schwarz preconditioners for nonlinear finite element simulations of laser beam welding processes}
\titlerunning{Monolithic Schwarz preconditioners in laser beam welding}  % abbreviated title (for running head)
\author{Tommaso Bevilacqua\inst{1} \and Axel Klawonn\inst{1,2} \and Martin Lanser\inst{1,2} \and Adam Wasiak\inst{1}}
%
%\authorrunning{Tommaso Bevilacqua et al.} % abbreviated author list (for running head)
%\tocauthor{Axel Klawonn, Martin Lanser, Adam Wasiak}
%
\institute{$^1$Department Mathematik/Informatik, Universit\"at zu K\"oln \\ Weyertal 86-90, 50931 K\"oln \\
$^2$Center for Data and Simulation Science,  Universit\"at zu K\"oln, Albertus-Magnus-Platz, 50923 K\"oln, Germany\\
\email{[tommaso.bevilacqua,axel.klawonn,martin.lanser,adam.wasiak]@uni-koeln.de}}

\maketitle
% typeset the title of the contribution
\begin{abstract}
Highly resolved finite element simulations of a laser beam welding process are considered. The thermomechanical behavior of this process is modeled with a set of thermoelasticity equations resulting in a nonlinear, nonsymmetric saddle point system. Newton's method is used to solve the nonlinearity and suitable domain decomposition preconditioners are applied to accelerate the convergence of the iterative method used to solve all linearized systems. Since a onelevel Schwarz preconditioner is in general not scalable, a second level has to be added. Therefore, the construction and numerical analysis of a monolithic, twolevel overlapping Schwarz approach with the GDSW (Generalized Dryja-Smith-Widlund) coarse space and an economic variant thereof are presented here.
\end{abstract}

\section{Introduction}
Laser beam welding is a noncontact joining process that has become considerably more important in the course of the increasing degree of automation in industrial production due to the fact that it requires short cycle times and small heat-affected zones. However, due to the high cooling rate inherent in the process that can lead to a residual melt oversaturated with certain alloy elements, solidification cracks can occur. A quantitative understanding of the development mechanisms of these fractures and their correlation to process parameters is thus indispensable for the development and improvement of this process. In particular, this is the focus of the DFG research group 5134 ''Solidification Cracks During Laser Beam Welding -- High Performance Computing for High Performance Processing'' \footnote[3]{https://www.for5134.science/en/}.

\begin{figure}[b!]
    \includegraphics[width=0.495\textwidth,trim={10cm 6cm 9cm 8cm},clip]{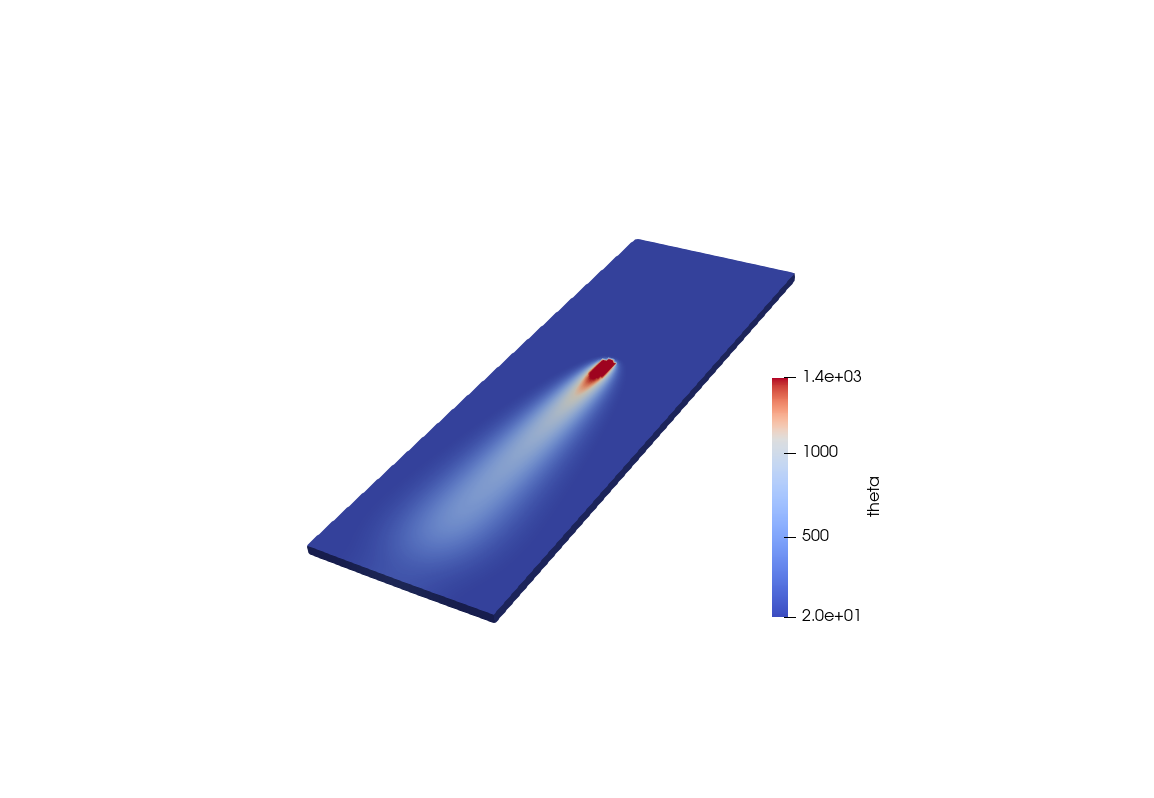}
    \includegraphics[width=0.495\textwidth,trim={10cm 6cm 9cm 8cm},clip]{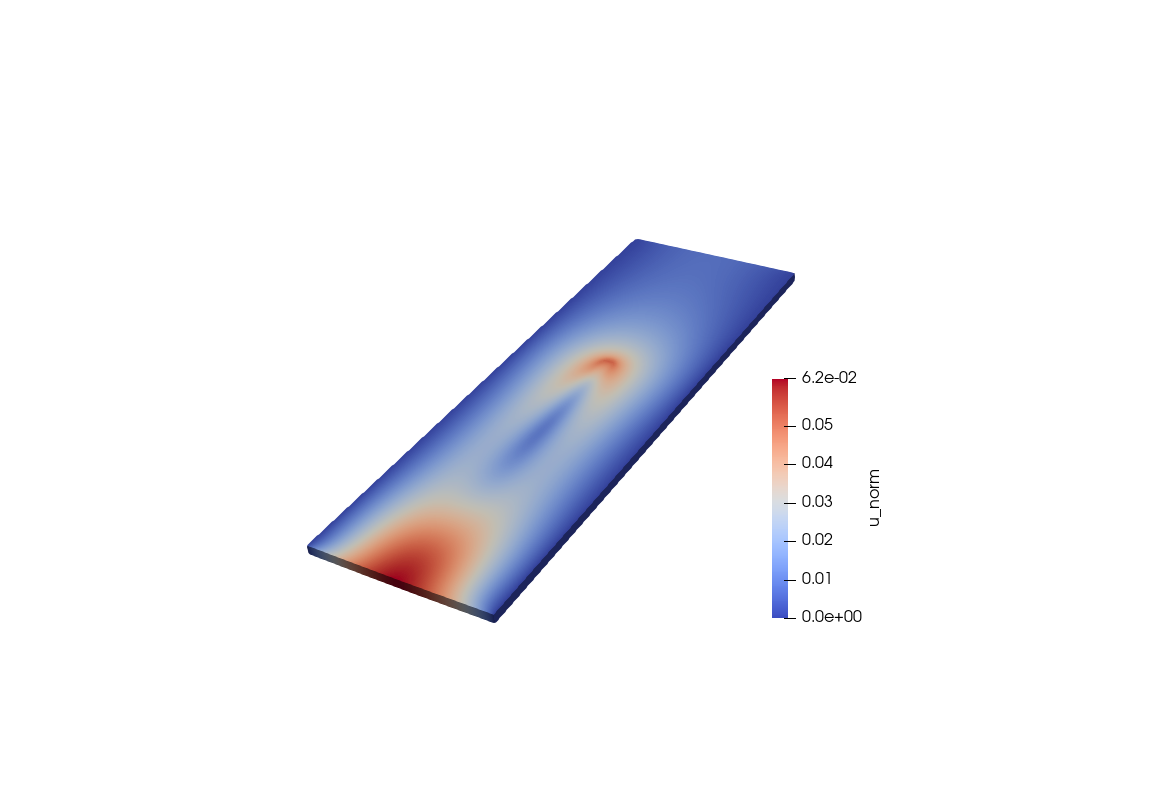}
{
  \caption{Example of temperature field (left) and norm of the displacements (right) computed after 100 time step iterations.}\label{displ}
}

\end{figure}

The thermomechanical properties of laser beam welding can be represented by a coupled system of time dependent thermoelasticity equations \cite{simo1992associative}. Due to the nature of the problem, the linear system that arises from the finite element discretization in space, the Backward Euler discretization in time, and the linearization by Newton's method is ill-conditioned. The scope of this work is to provide a first choice of a preconditioner to accelerate the convergence of the GMRES (Generalized Minimal Residual) method. Since we are predominantly interested in analyzing the performance of the preconditioners, we consider a  model problem where the action of the laser beam is simplified. Therefore, we %load 
obtain the geometry of the melting pool, that is, the area of material which is melted by the laser, from experimental data \cite{bakir2018numerical} and enforce the melting temperature of the specific material in all degrees of freedom within this area using  Dirichlet boundary conditions. \newline
\indent We note that the laser beam and thus, the melting pool geometry and corresponding Dirichlet boundaries move over time with a problem dependent speed. We have chosen this approach with future work in mind, where a realistic heat source model for the laser will be used, in order to obtain more accurate simulations \cite{hartwig2024physically,hartwig2023numerical,hartwig2023volumetric}. For completeness, \cref{displ} gives a rough overview of a typical model geometry which is welded by a laser and \cref{meltpool} shows the geometry of the melting pool we use. \newline 
\indent In this work we make use of a monolithic overlapping Schwarz preconditioner for saddle point problems. 
In order to obtain numerical scalability, a coarse level needs to be introduced. 
 Lagrangian coarse spaces that operate on the saddle point problem have been introduced by Klawonn and Pavarino in \cite{klawonn1998overlapping,klawonn2000comparison}. To obtain more flexibility in the construction of the coarse space we use variants of the Generalized Dryja-Smith-Widlund (GDSW) coarse spaces, in which basis functions are defined using problem dependent extension operators; in case of symmetric positive definite matrices these extensions are discrete harmonic. This coarse space originally introduced in \cite{dohrmann2008domain,dohrmann2008family} has been also extended to saddle point problems arising from the discretization of incompressible fluid flow problems \cite{heinlein2019monolithic}. 
Since the linearized systems of fluid flow and thermo-elastic problems have a similar block structure with similar coupling blocks, we apply and slightly customize the preconditioners from \cite{heinlein2019monolithic} here and numerically analyze them for our thermoelastic laser beam welding problems. We additionally propose an economic variant that builds the coarse space by solving uncoupled problems instead of the full saddle point problem. \newline 
\indent The remainder of the paper is organized as follows: In \cref{thermo-eq} we present the thermoelasticity equations, in \cref{FEMdisc} we introduce the discretization, in \cref{SchPrec} we present the preconditioner, and finally in \cref{NumRes} we show some numerical results.
\section{Thermoelasticity equations}\label{thermo-eq}
We describe the thermomechanical behavior of a plate of metal welded by a laser with the set of thermoelasticity equations
\begin{equation}
  \label{eq:1}
  \begin{split}
 {\rm div} ({\sigma(u,\theta)}) + \rho \, \mathbf{b} &= 0, \\ 
 {\rm div} (\mathbf{q})+\gamma \, {\rm div} (\dot{u}) \, \theta +c_\rho\, \dot{\theta} &= \rho \, r,
 \end{split}
\end{equation}
where $u(x,t)$ and $\theta(x,t)$ represent the displacement vector and the absolute temperature respectively. The other variables are: $r$ the heat source, $\rho$ the mass density, $c_\rho$ the heat capacity, $\mathbf{b}$ the volume acceleration, $\mathbf{q}$ the heat flux vector, and $\gamma = \alpha_T(3\lambda + 2\mu) = 3\alpha_T\kappa$ the stress temperature modulus. Moreover, $\kappa = \lambda + 2\mu/3$ denotes the bulk modulus, $\alpha_T$ the coefficient of linear thermal expansion, and $\lambda, \mu$ the Lam\'e parameters. 
The stress in tensor notation is defined as
\begin{equation*}
\mathbf{\sigma} := \mathbb{C} : \varepsilon- 3\alpha_T \kappa (\theta -\theta_0) \bm{1}
\end{equation*}
with $\mathbb{C} = \kappa \bm{1} \otimes \bm{1} + 2\mu \mathbb{P}$ and the fourth-order deviatoric projection tensor $\mathbb{P} = \mathbb{I} - \frac{1}{3} \bm{1} \otimes \bm{1}$. \\
The heat flux vector is described by Fourier's law $\mathbf{q} = -\bm{\lambda} \cdot \nabla \theta$, where $\bm{\lambda}$ is the identity matrix multiplied by the heat conduction coefficient. For the sake of simplicity, we assume to not have any volume acceleration, that is, $\mathbf{b}=0$ and no heat source, therefore $r=0$. The only heat source is the laser  which is modeled by a volume Dirichlet boundary condition that represents the pool of melted metal that moves along the plate. This geometry is %loaded 
obtained from experimental data %that give us a 
which yields a triangulated surface \cite{bakir2018numerical}. 
Thus, we obtain the following simplified form of~\cref{eq:1}
\begin{align*}
 {\rm div} ({\sigma(u,\theta)}) &= 0, \\
 {\rm div} (\mathbf{q})+\gamma \, {\rm div} (\dot{u}) \, \theta +c_\rho \, \dot{\theta} &= 0.
\end{align*}
The weak formulation of these equations is
\begin{align*}
 \int_\Omega \sigma (\mathbf{u})\, :\, \varepsilon (\mathbf{v}) \diff x &= 0 \qquad \forall \mathbf{v} \in [H^1(\Omega)]^3, \\
\int_\Omega \mathbf{q} \cdot \nabla q \diff x + \int_\Omega \{- 3\alpha_T \, \kappa \,{\rm tr[\dot{\varepsilon}] \, \theta - c_\rho \dot{\theta}} \}q \diff x &= 0 \qquad \,\forall q \in H^1(\Omega),
\end{align*}
with Dirichlet boundary conditions $\mathbf{u} = 0$ on $\partial \Omega_{D,u} \subset \partial \Omega$ and $\theta = \theta_{l} \gg~\!\!0$ on $\partial \Omega_{D,\theta} \subset \partial \Omega$. This is a nonlinear and nonsymmetric saddle point system for which the theoretical framework is still an object of studies, for example, an investigation of pairs of inf-sup-stable finite elements. A theoretical investigation of this problem is beyond the scope of this work. 

\section{Finite element discretization}\label{FEMdisc} 
Applying the Newton-Raphson method to the nonlinear boundary value problem yields, in each nonlinear iteration, a linear problem which spatially is discretized by using a mixed finite element method; see \cite{hartwig2024physically} for further details on the linearization. Next, we consider the mixed finite element discretization. 
%\mycomment{AK: Reformulated. OK this way? Check order linearization and time stepping.\\
%TB: Yes, thanks! I added the definition of the time step since it was missing.}
%We solve the system of equations resulting from the calculation of the boundary value problem by applying the Newton-Raphson iteration method (see \cite{hartwig2024physically} for further details) and perform the spatial discretization of the saddle point problem using the mixed finite element method. 

Let $\tau_h$ be a uniform mesh of $N_e$ hexahedral elements $\Omega_e$ of $\Omega$ with characteristic mesh size $h$. We introduce the conforming discrete piecewise linear displacement and temperature spaces
\begin{align*}
V^h &= V^h(\Omega) = \{\mathbf{u}\in [\mathcal{C}^0(\Bar{\Omega}) \cap H^1(\Omega)]^3: \mathbf{u}|_T \in Q_1 \ \forall T \in \tau_h\}, \\
Q^h &= Q^h(\Omega) = \{\theta \in \mathcal{C}^0(\Bar{\Omega}) \cap H^1(\Omega): \theta|_T \in Q_1 \ \forall T \in \tau_h\}
\end{align*}
respectively, of Q1-Q1 mixed finite elements. Here, $\mathcal{C}^0(\Bar{\Omega})$ denotes the space of continuous functions on $\Bar{\Omega}$ and $H^1(\Omega)$ the usual Sobolev space. Since we do not have a stable theoretical framework, we can not ensure that this choice of finite elements satisfies an inf-sup condition \cite{boffi2013mixed}. The time discretization is done by using a Backward Euler method with time step $\Delta t$. Therefore, in each time step, we have to solve a linearized system until the global absolute residual of Newton's method falls below a fixed tolerance. The resulting discrete saddle point problem has the form
\begin{equation*}
    K\varDelta \mathbf{d} = \begin{bmatrix}
        K_{uu} & K_{u\theta} \\
        K_{\theta u} & K_{\theta \theta}
    \end{bmatrix} \begin{bmatrix}
        \varDelta \mathbf{d}_u \\ \varDelta \mathbf{d}_\theta 
    \end{bmatrix} = \begin{bmatrix}
        R_u \\ R_\theta 
    \end{bmatrix} = R, 
\end{equation*}
where $\varDelta \mathbf{d}_u$ and $\varDelta \mathbf{d}_\theta$ represent the Newton update for the displacement and temperature, and $R_u$ and $R_\theta$ the vectors of the residual, respectively. The block matrices are obtained by finite element assembly and we obtain
\begin{equation*}
\begin{split}
& K_{uu} = \mathbb{A}^{N_e}_{e=1} \bigg[ \int_{\Omega_e} \mathbf{B}_u^T \, \mathbb{C}\, \mathbf{B}_u \diff x \bigg], \qquad K_{u\theta} = \mathbb{A}^{N_e}_{e=1} \bigg[ -\int_{\Omega_e} \mathbf{B}_u^T \, (3\alpha_T \kappa \mathbf{1})\, {\rm \mathbf{N}}_\theta \diff x \bigg],\\
& \hspace{20mm} K_{\theta u} = \mathbb{A}^{N_e}_{e=1} \bigg[-\frac{1}{\varDelta t} \int_{\Omega_e} \theta \, {\rm \mathbf{N}}_\theta^T \, (3\alpha_T \kappa \mathbf{1}) \, \mathbf{B}_u \diff x \bigg], \\ 
& K_{\theta\theta} = \mathbb{A}^{N_e}_{e=1} \bigg[ -\int_{\Omega_e} \mathbf{B}_\theta^T \, \bm{\lambda} \, \mathbf{B}_\theta \diff x \, - \int_{\Omega_e} {\rm \mathbf{N}}_\theta^T \, {\rm tr}[\dot{\varepsilon}] \, {\rm \mathbf{N}}_\theta \diff x \,- \frac{1}{\varDelta t}\int_{\Omega_e} {\rm \mathbf{N}}_\theta^T \, c_\rho \, {\rm \mathbf{N}}_\theta \diff x\bigg],
\end{split}
\end{equation*}
where $\mathbf{N}_u$ and ${\rm \mathbf{N}}_\theta$ are the finite element nodal basis functions for displacement and temperature and $\mathbf{B}_u$ and $\mathbf{B}_\theta$ denote their derivatives in tensor notation.
The nature of this linear system requires suitable preconditioners to accelerate the convergence of an iterative solver, as, in our case, the GMRES method. Our choice is to use one- and two-level monolithic Schwarz preconditioners and a computationally cheaper variant that uses a block diagonal approach to build the coarse level.
\section{Schwarz preconditioners}\label{SchPrec}
\subsection{Monolithic two-level Schwarz preconditioner}
We consider a monolithic overlapping Schwarz domain decomposition preconditioner of GDSW (Generalized Dryja-Smith-Widlund) type for the saddle point problem. This is a two-level Schwarz preconditioner where the coarse space and the local problems have the same block structure as the original one. We refer the interested reader to \cite{dohrmann2008domain,dohrmann2008family} for a complete definition of the GDSW preconditioner for elliptic problems and to  \cite{heinlein2019monolithic} for saddle point problems resulting from Stokes' equations.

We introduce $\{ \Omega_i\}^N_{i=1}$, a nonoverlapping domain decomposition of the domain $\Omega$ into $N$ subdomains of characteristic diameter $H$, and $\{ \Omega'_i\}^N_{i=1}$, the corresponding overlapping domain decomposition with overlap $k$, generated by geometric or graph techniques. We define the interface $\Gamma$ of the nonoverlapping domain decomposition
$\Gamma=\{ x\in (\Omega_i \cap \Omega_j ) \setminus \partial \Omega_D : i \neq j, 1\leq i, j\leq N \},$ that is divided into $M$ connected components, specifically into vertices, edges, and faces, which belongs to the same set of subdomains. These edges, faces, and vertices are used later on to define the basis functions of the GDSW coarse space. In general, the definition of an interface allows for the partitioning of the degrees of freedom into interface ($\Gamma$) and interior ($I$) ones and, after a reordering, an equivalent partitioning of $K$, that is,  
\begin{equation*}
K:= \begin{bmatrix}
        K_{II} & K_{I\Gamma} \\
        K_{\Gamma I} & K_{\Gamma \Gamma}
    \end{bmatrix}.
\end{equation*}
We further introduce the operators $R_{u,i} : V^h \rightarrow V^h_i$ and $R_{\theta,i} : Q^h \rightarrow Q^h_i$ for $i = 1,...,N$. These are restrictions from the global finite element spaces to the local finite element spaces defined on the overlapping subdomains. The monolithic restriction operators are then denoted by $R_i : V^h \times Q^h \rightarrow V^h_i \times Q^h_i$ for $i=1,...,N$, with
\begin{equation*}
R_i:= \begin{bmatrix}
        R_{u,i} & 0 \\
        0 & R_{\theta,i}
    \end{bmatrix}.
\end{equation*} 
The transposed operators $R^T_i$, $R^T_{u,i}$, and $R^T_{\theta,i}$ are the corresponding prolongation operators.
Now, the additive monolithic two-level GDSW preconditioner can be written as
\begin{equation*}
\hat{B}^{-1}_\text{GDSW} = \Phi K^{-1}_0 \Phi^T + \sum_{i=1}^N R_i^T \, K_i^{-1} \, R_i,
\end{equation*}
where the local stiffness matrices $K_i$ are extracted from $K$ by
\begin{equation*}
K_i = R_i K R_i^T , \qquad {\rm for}\,\, i = 1, . . . , N,
\end{equation*} 
and the coarse operator by
\begin{equation*}
K_0=\Phi^T K \Phi.
\end{equation*}
Here, $\Phi$ is the matrix that contains the coarse space basis functions $\phi^i$ for $i = 1,...,N_\phi$ column-wise. These functions are constructed as the discrete saddle point extensions of functions $\phi^i_\Gamma$ defined on the interface by solving the linear system
\begin{equation}\label{interfprb}
\begin{bmatrix}
        K_{II} & K_{I \Gamma} \\
        0 & I
    \end{bmatrix} \begin{bmatrix}
        \phi^{i}_I \\ \phi^i_{\Gamma}
    \end{bmatrix} = \begin{bmatrix}
        0 \\ \phi^i_{\Gamma}
    \end{bmatrix},
\end{equation}
where we decompose the interface values of the coarse basis into displacements ($u_0$) and temperature ($\theta_0$) as $\phi^i_{\Gamma} = [ \,\phi^{i\,T}_{\Gamma,u_0} \ \phi^{i\,T}_{\Gamma,\theta_0}\,]^T$.
The entries of $\phi^i_{\Gamma}$ are usually set to 1 or 0 and generally represent the restriction of the null space of the operators $K_{uu}$ and $K_{u \theta}$ on the interface components (edges, vertices, and faces) $\Gamma_j$ with $j=1,...,M$. 

In three dimensions the null space of the $K_{uu}$ block is given by the six rigid body modes (three translations and three rotations) defined on the interface components, but here we only consider the translations. For the $K_{u \theta}$ operator, these consist of the functions that have constant temperature on $\Omega$, so we restrict the constant function 1 to vertices, edges, and faces.

Specifically, each node of our mesh can be associated to a vector with four components, where the first three entries represent the displacements and the last one the temperature. We introduce the vectors
\begin{equation*}
    v_{u,x}=\begin{bmatrix}1 \ 0 \ 0 \ 0\end{bmatrix}^T, \quad
    v_{u,y}=\begin{bmatrix}0 \ 1 \ 0 \ 0\end{bmatrix}^T, \quad
    v_{u,z}=\begin{bmatrix}0 \ 0 \ 1 \ 0\end{bmatrix}^T, \quad
    v_{\theta}=\begin{bmatrix}0 \ 0 \ 0 \ 1\end{bmatrix}^T
\end{equation*}
where the first three are the translations into $x$-, $y$-, and $z$- directions and the last one is the restriction of the constant function on the node. For each vertex, edge, and face $\Gamma_j,\;j=1,...,M,$ of the interface $\Gamma$, we then construct four basis functions  $\phi_{\Gamma_j}^k$, $k=1,...,4$, by filling the degrees of freedom associated to the nodes of $\Gamma_j$ with $v_{u,x}$ for $j=1$, $v_{u,y}$ for $j=2$, etc. For the nodes in $\Gamma \setminus \Gamma_j$, $\Phi_{\Gamma_j}^k$ is set to zero. Finally, we collect all $\Phi_{\Gamma_j}^k,\; j=1,...,M, \; k=1,...,4$ as the columns of $\Phi_\Gamma$ and thus obtain $4M$ coarse basis functions and columns in $\Phi$.

After having solved \cref{interfprb} for each basis function of the coarse space, the matrix $\Phi$ has the following block representation
\begin{equation*}
\Phi = \begin{bmatrix}
        \Phi_{u,u_0} & \Phi_{u,\theta_0} \\
        \Phi_{\theta_0, u} & \Phi_{\theta,\theta_0}
    \end{bmatrix}.
\end{equation*}
With this construction, the presence of nonzero blocks on the off-diagonal can deteriorate the numerical scalability, especially if nonzero Dirichlet boundary conditions for the temperature are present (cf. \cite{heinlein2019monolithic}), as we have with the laser beam. Therefore, these blocks need to be removed, obtaining 
\begin{equation*}
\Phi = \begin{bmatrix}
        \Phi_{u,u_0} & 0 \\
        0 & \Phi_{\theta,\theta_0}
    \end{bmatrix}.
\end{equation*}

\subsection{Economic variant of the two-level Schwarz preconditioner}
We present here an economic variant of the construction of the GDSW coarse space. In order to construct the first level of the preconditioner we solve the saddle point system on each subdomain, whereas the construction of the coarse basis function is obtained by solving two uncoupled problems. This can be a valuable option in future works to save computational time.
The resulting preconditioner can be written as
\begin{equation}
\hat{B}^{-1}_\text{EGDSW} = \hat{\Phi} \hat{K}^{-1}_0 \, \hat{\Phi}^T + \sum_{i=1}^N R_i^T K_i^{-1}  R_i,
\end{equation}
where the subdomain local problems are solved as before, while the coarse operator is defined as
\begin{equation*}
\hat{K}_0=\hat{\Phi}^T \hat{K} \hat{\Phi},\qquad\text{where}\qquad\hat{K} = \begin{bmatrix}
        K_{uu} & 0 \\
        0 & K_{\theta \theta}
    \end{bmatrix}.
\end{equation*}
The choice of the interface basis functions and the extension to the subdomain interiors is performed as before by solving \eqref{interfprb}, but now $K_{II}$ and $K_{I\Gamma}$ do not contain coupling blocks, or, in other words, we remove the off diagonal blocks before performing the extension and thus use a block diagonal approach here. The final prolongation matrix $\hat{\Phi}$ is consequently of block diagonal form
\begin{equation*}
\hat{\Phi} = \begin{bmatrix}
        \hat{\Phi}_{u,u_0} & 0 \\
        0 & \hat{\Phi}_{\theta,\theta_0}
    \end{bmatrix}.
\end{equation*}

\section{Numerical results}\label{NumRes}
\begin{figure}[t!]
\begin{floatrow}
\capbtabbox{%
\caption{Parameters of the material austenitic chrome nickel steel(1.4301) at $20^{\circ}$C \cite{richter2010physikalischen}.}\label{TabelCoefficient}
  \begin{tabular}{r@{\quad}rl}
\hline
\multicolumn{1}{l}{\rule{0pt}{12pt}
Parameter}&\multicolumn{1}{l}{Values}\\[2pt]
\hline\rule{0pt}{12pt}
$E$  &  198000 ${\rm N/mm}^2$ \\
$\nu$ & 0.276 ${\rm N/mm}^2$ \\
$\alpha_T$&   $1.6 \cdot 10^{-5} \,{\rm K}^{-1}$ \\
$\rho$  &  $7919 \,{\rm kg/m}^3$ \\
$c_{\rho}$  & $468 \, {\rm J/kg}\cdot {\rm K}$ \\
$\bm{\lambda}$ & $(14.4 \, {\rm W}/({\rm m}\cdot {\rm K})) \cdot \bm{I}$ \\[2pt]
\hline
\end{tabular}
}{%
  
}\quad
\ffigbox{%
    \includegraphics[width=0.27\textwidth]{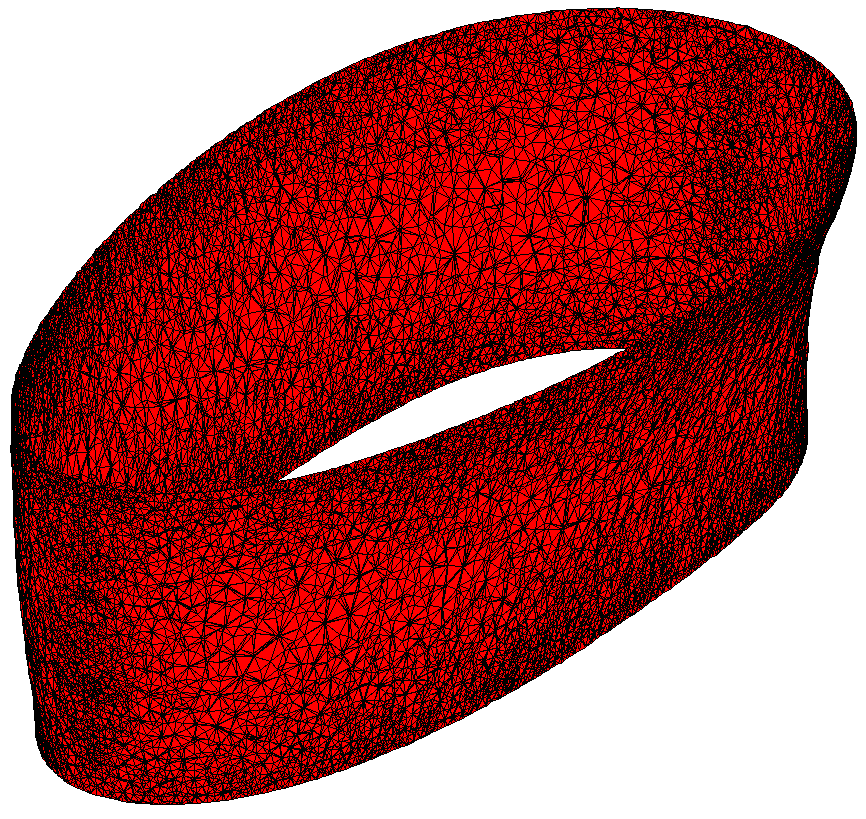}
}{%
  \caption{Triangulated surface representing the geometry of the melting pool \cite{bakir2018numerical}.}\label{meltpool}
}
\end{floatrow}
\end{figure}

The numerical simulations are performed by an internal C/C++ software with an MPI/OpenMP~\cite{chandra2001parallel} hybrid parallel implementation exploiting the PETSc library 
\cite{balay2019petsc}. We equipped our software with an interface to the FEAP library \cite{taylor2014feap}, to exploit the implemented finite elements and material models from FEAP. The local and coarse problems are solved exactly using MUMPS \cite{amestoy2000mumps}.
 We solve a thermomechanical problem on a domain $\Omega$ = $60 \, {\rm mm} \times 20\, {\rm mm} \times 1 \,{\rm mm}$, {where} $(0,0,0)$ is one of the corners enforcing Dirichlet boundary conditions for the displacements on the faces $y = 0\,{\rm mm}$ and $y = 20\,{\rm mm}$. Our configuration models a laser moving along a plate of an austenitic chrome nickel steel; 
see \cref{TabelCoefficient} for the chosen material parameters. We note that for this type of simulation we are using constant material parameters. An introduction of temperature dependent parameters can slightly, but not significantly, affect the performance of the GDSW preconditioners. 
 The action of the laser is implemented as a volume Dirichlet boundary condition 
 which represents the boundary of the region of the metal that has been melted (melting pool)  as described before; see \cref{meltpool} for the geometry of the melting pool. Therefore, we load a surface model obtained from experimental data that models the region of the melting pool and we move this geometry along the $x$-direction during the time iterations at a certain speed. The initial temperature of the metal is set to $\theta_0 = 20^\circ $C everywhere and the temperature inside of the melting pool is set, as a gradient, to $\theta = 1480^\circ $C.

\begin{table}[b!]
\caption{Weak and strong scalability test for a fully coupled thermoelasticity problem modeling a laser beam welding process with Q1-Q1 elements. For the weak scaling experiments the number of elements for each subdomain is fixed to $10 \times 5 \times 10$. The mark "x" means that the method does not converge within the stopping criteria.}\label{ThreeSchwarzW}
\begin{center}
\begin{tabular}{crrrrrrrr}
\hline
& \multirow{2}{*}{nProcs} & \multirow{2}{*}{nDoFs} & \multicolumn{2}{c}{one-level} & \multicolumn{2}{c}{GDSW} & \multicolumn{2}{c}{EGDSW} \\[-1.5mm]
& & &  \multicolumn{1}{c}{${\rm it}_\text{Avg}$} & \multicolumn{1}{c}{${\rm it}_{N}$} & \multicolumn{1}{c}{${\rm it}_\text{Avg}$} & \multicolumn{1}{c}{${\rm it}_{N}$} &\multicolumn{1}{c}{${\rm it}_\text{Avg}$} & \multicolumn{1}{c}{${\rm it}_{N}$}  \\
\hline\hline 
&$4\times4\times1$   & \numgru{37884}   & 14 & 4 & 11 & 4 & 11 & 4 \\
\bf{Weak} &$8\times8\times1$   & \numgru{146124}  & 26 & 4 & 16 & 4 & 17 & 4 \\
\bf{Scaling} &$16\times16\times1$  & \numgru{573804} & 284 & 4 & 89 & 4 & 89 & 4 \\
&$32\times32\times1$ & \numgru{2273964} & x & x & 120 & 4 & 122 & 4	\\[-0.3mm]
\hline
\hline
&$4\times4\times1$                & \numgru{2273964} & 140 & 4 & 102 & 4 & 104 & 4 \\
\bf{Strong} &$8\times8\times1$    & \numgru{2273964} & 285 & 4 & 113 & 4 & 115 & 4 \\
\bf{Scaling} &$16\times16\times1$ & \numgru{2273964} & 466 & 4 & 110 & 4 & 116 & 4 \\
&$32\times32\times1$              & \numgru{2273964} & x & x & 120 & 4 & 122 & 4	\\[-0.3mm]
\hline
\end{tabular}
\end{center}

\end{table}
Since our results are preliminary and a detailed performance optimization of the parallel PETSc implementation is still in progress, we limit ourselves to consider the iteration counts averaged over the first three time iterations with a step size of $\varDelta t = 0.05$ and omit any timings or parallel execution times for now. Still, these will be of interest in the future. \\
Let us now introduce the legend that we will use in all tables: nProcs = number of cores (equal to the number of subdomains), nDoFs = number of degrees of freedom of the system, ${\rm nDofs}_\Gamma$ = number of degrees of freedom in the coarse space, ${\rm it}_\text{Avg}$ = rounded average number of GMRES iterations to solve each linear system, ${\rm it}_{N}$ = number of Newton iterations, and ${\rm it}_\text{Tot}$ = total number of GMRES iterations. The stopping criterion for  GMRES is either a residual reduction of $10^{-6}$ of the unpreconditioned residual or 1000 iterations, while an absolute residual of $10^{-8}$ or a maximum of 10 ${\rm it}_\text{N}$ is used as stopping criterion for Newton's method. 

\textbf{Strong and Weak Scalability}.
In \cref{ThreeSchwarzW} we compare the iteration counts for three different choices of Schwarz preconditioners. As expected, the one-level Schwarz method does not scale, while the two GDSW approaches show a better behavior. We note that the number of iterations here is also affected by a jump in the temperature of $10^3$ introduced by the laser and the high temperature in the melting pool. 

\begin{table}[t!]
\caption{Weak scalability test for a fully coupled thermoelasticity problem with Q1-Q1 elements. For the configuration with only one subdomain in $z$-direction the local number of elements is fixed to $5 \times 5 \times 10$. For the one with two subdomains in $z$-direction it is fixed to $10 \times 5 \times 5$.}\label{WeakCoarse}
\begin{center}
\begin{tabular}{rrrrrrrrr}
\hline
 \multirow{2}{*}{nProcs} & \multirow{2}{*}{nDoFs} & \multirow{2}{*}{${\rm nDoFs}_\Gamma$} &\multicolumn{3}{c}{GDSW} & \multicolumn{3}{c}{EGDSW} \\[-1.5mm]
& & & \multicolumn{1}{l}{${\rm it}_\text{Avg}$} & \multicolumn{1}{l}{${\rm it}_{N}$} & \multicolumn{1}{l}{${\rm it}_\text{Tot}$} & \multicolumn{1}{l}{$it_\text{Avg}$} & \multicolumn{1}{l}{${\rm it}_{N}$} & \multicolumn{1}{l}{${\rm it}_\text{Tot}$} \\
\hline\hline
$32\times8\times1$   & \numgru{290444} & \numgru{2756} & 14 & 4 & 57 & 14 & 4 & 57 \\
$64\times16\times1$  & \numgru{1140444} & \numgru{11652} & 82 & 4 & 328 & 82 & 4 & 328 \\
$128\times32\times1$ & \numgru{4540844} & \numgru{47876} & 121 & 4 & 486 & 119 & 4 & 476	\\[-0.3mm]
\hline
\hline
$16\times8\times2$   & \numgru{290444} & \numgru{4556} & 28 & 4 & 110 & 30 & 4 & 118 \\
$32\times16\times2$  & \numgru{1140444} & \numgru{19430} & 127 & 4 & 510 & 145 & 4 & 582 \\
$64\times32\times2$ & \numgru{4540844} & \numgru{79628} & 110 & 4 & 440 & 111 & 4 & 445	\\[-0.3mm]
\hline
\end{tabular}
\end{center}
\end{table}

\textbf{Weak Scalability with different partitionings}.
In \cref{WeakCoarse} we reproduced a larger weak scaling test for two different configurations of subdomain partitions. Both coarse spaces seem to perform well with, as expected, generally low iteration counts for the monolithic GDSW coarse space. Since the difference is not so evident we will also apply and analyze the economic variant of the coarse space in future studies, in particular with regard to computing time. To summarize, both monolithic two-level Schwarz preconditioners perform well with respect to iteration counts for the laser beam welding processes even though the high jump in the temperature, the thin plate geometry, and the lack of theory introduce some severe difficulties.

\paragraph{Acknowledgement}
This project has received funding from the Deutsche Forschungsgemeinschaft (DFG) as part of the Forschungsgruppe (Research Unit) 5134 ``Solidification Cracks During Laser Beam Welding -- High Performance Computing for High Performance Processing''.
The authors gratefully acknowledge the scientific support and HPC resources provided by the Erlangen National High Performance Computing Center (NHR@FAU) of the Friedrich-Alexander-Universit\"at Erlangen-N\"urnberg (FAU) under the NHR project b144dc. NHR funding is provided by federal and Bavarian state authorities. NHR@FAU hardware is partially funded by the German Research Foundation (DFG) - 440719683. We would like to thank our project partners from Forschungsgruppe 5134 L.\ Scheunemann, J. Schr\"oder, and P. Hartwig for support with the interface to FEAP and providing the thermoelasticity formulation and all parameters in FEAP and also M. Rethmeier and A. Gumenyuk for providing the surface data of the melting pool.

\bibliographystyle{siam}
\bibliography{literature}

\end{document}